\documentclass[11pt, reqno]{amsart}
\usepackage[margin=1.2in]{geometry}                % See geometry.pdf to learn the layout options. There are lots.
\usepackage{graphicx}
\usepackage{amssymb}
\usepackage{epstopdf}
\usepackage{tikz-cd}
\usepackage{forest}
\usepackage[colorlinks=true, pdfstartview=FitV, linkcolor=blue, citecolor=blue, urlcolor=blue]{hyperref}
 \usepackage{tikz}
\usepackage[listings]{tcolorbox}
\usepackage{caption}
\usepackage{subcaption}
\usepackage{mathrsfs}
\usepackage{paralist}
\usepackage{verbatim}
\usepackage{animate}
\usepackage[smalltableaux]{ytableau}
\usepackage[nameinlink]{cleveref}
\usepackage{accents}
\crefname{defn}{definition}{definitions}
\Crefname{def}{Definition}{Definitions}
\crefname{thm}{theorem}{theorems}
\Crefname{thm}{Theorem}{Theorems}

\usepackage{color}

\newcommand{\N}{{\mathbb N}}

\newcommand{\Q}{{\mathbb Q}}

\newcommand{\M}{{\mathcal M}}
\renewcommand{\O}{{\mathcal O}}
\renewcommand{\P}{{\mathcal P}}

\newcommand{\cP}{{\mathcal P}}

\newcommand{\dSdw}{\Delta}
\newcommand{\uSdw}{\rotatebox[origin=c]{180}{$\Delta$}}

\def\Seg{\operatorname{Seg}}

\def\Span{\operatorname{Span}}

\theoremstyle{plain} %% This is the default, anyway
\newtheorem{thm}{Theorem}[section]

\newtheorem*{introthm*}{Theorem}

\newtheorem*{oprobl*}{Open Problem}

\theoremstyle{definition}
\newtheorem{defn}[thm]{Definition}

\newtheorem{ex}[thm]{Example}

\theoremstyle{remark}

\numberwithin{equation}{section}  %% equation numbering

\title{The MacaulayPosets package for Macaulay2}

\author{Penelope Beall}
\address{University of California--Davis, 1 Sields Avenue, Davis, CA 95616, United States}
\email{pbeall@ucdavis.edu}

\author{Nikola Kuzmanovski}
\address{University of Notre Dame, 138 Hayes-Healy, Notre Dame, IN, 46556-4618, United States}
\email{nkuzmano@nd.edu}

\author{Yu Olivier Li}
\address{University of St Andrews, Mathematical Institute, St Andrews KY16 9SS, United Kingdom}
\email{yuolivierli@gmail.com}

\author{Alexandra Seceleanu}
\address{University of Nebraska--Lincoln, 203 Avery Hall, Lincoln, NE 68588, United States}
\email{aseceleanu@unl.edu}

\thanks{We acknowledge the support of NSF DMS--2341670
for the Polymath Jr.\,program, which facilitated our collaboration. 
Seceleanu was partially supported by NSF grant DMS--2401482.}

\begin{document}

\begin{abstract}
    We introduce the \texttt{Macaulay2} package \texttt{MacaulayPosets}. This package utilized the poset data type introduced in the \texttt{Posets} package and offers functionality for studying the Macaulay property for posets, particularly those which arise as monomial posets of commutative rings. 
\end{abstract}

\maketitle
\vspace{-2em}
\section{Introduction}
Partially ordered sets, or posets, are fundamental combinatorial structures widely used in contemporary research to model discrete systems endowed with an order relation. Macaulay2 provides the \texttt{Posets} package \cite{PosetsJSAG}, developed by David Cook II, Sonja Mapes, and Gwyneth Whieldon, which introduces a dedicated data structure along with essential methods for working with posets. %This package includes tools for enumerating commonly studied classes of posets, performing various operations, and computing key invariants associated with posets.

The development of the \texttt{MacaulayPosets} package is motivated by the study of posets in the context of commutative algebra. This package builds upon the poset data type introduced in the \texttt{Posets} package, extending its functionality to explore the Macaulay property in posets, particularly those that arise as monomial posets of commutative rings. A Macaulay poset is characterized by a ranked structure and a total order that interacts harmoniously with the partial order, enabling the establishment of bounds on the sizes of subsets of a given rank within an order ideal. Prior to our package, there was a single software \cite{SoftwareTool} developed in Java by S.\,L.\, Bezrukov and A.\,Dissanayake, which could check whether a poset with at most $30$ elements is Macaulay. Our package surpasses this limitation (see \Cref{ex: 39 vertices}) and allows to check the Macaulay property for graded rings in addition to posets.

Macaulay posets provide a framework for understanding the interplay between algebraic invariants of graded rings, such as Hilbert functions, and their combinatorial properties. Their foundational application, introduced by F. S. Macaulay in \cite{Mac}, was to characterize all possible Hilbert functions of homogeneous ideals in polynomial rings. Today, Macaulay posets play a  role in both commutative algebra and enumerative combinatorics, making them a rich subject of research. In extremal combinatorics, they appear in the study of isoperimetric problems, which ask how large or small the set of neighbors of a finite set can be, given its size, within a discrete structure such as a graph. Significant contributions to these problems using Macaulay posets were made by  Katona \cite{Katona}, Kruskal \cite{Kruskal}, and Clements-Lindstr\"om \cite{CL}.  Surveys on Macaulay posets in combinatorics can be found in \cite{BL} and \cite[Chapter 8]{Engel}.
Subsequent to the combinatorial developments, an emerging theory of Macaulay rings -- rings whose monomial poset exhibits the Macaulay property -- has gained attention in commutative algebra. Works such as \cite{Mermin, MP, MP2, Kuz} aim to extablish the Macaulay property for certain classes of rings, primarily monomial or toric quotients of a polynomial ring, leveraging additional algebraic structure or ring operations such as tensor product.  %This growing connection between algebra and combinatorics highlights the significance of Macaulay posets as a valuable tool in contemporary mathematical research.

Our package was conceived as a companion to the paper \cite{Polymath paper}, which considers when the Macaulay property is preserved under further poset and ring operations. As such, the package implements four poset operations dubbed wedge product, closed product, fiber product, and connected sum, and two related ring operations termed fiber product and connected sum. These operations allow to easily construct complex posets, some of which have the Macaulay property and some which do not.

\smallskip

\paragraph{\bf Acknowledgements.} We thank Sergei L.\,Bezrukov for sharing code for the software \cite{SoftwareTool} with us. While our package does not rely on this code, it gave us our first exposure to the computational  aspects of Macaulay posets.

\section{Macaulay Posets and Macaulay Rings}

In this section we give an overview the Macaulay property as it pertains both to posets and to rings. 

Throughout $\mathbb{N}$ denotes the set of non-negative integers, including zero.

\subsection{Macaulay posets} 
A {\bf partial order} on a set $\P$ is a relation on $\P$ that is antisymmetric, reflexive, and transitive. A {\bf partially ordered set} or {\bf poset} is a set $\P$ along with a partial order on $\P$. A partial order $\leq$ on $\P$ is a {\bf total order} if and only if for all $p,q\in\P$, we have $p\leq q$ or $q\leq p$.

Suppose $\leq$ is a partial order on $\P$. Suppose $p,q\in\P$. We say $q$ {\bf covers} $p$ if and only if $p<q$, and $p<r<q$ for no $r\in\P$.

\begin{defn}\label{def: shadows}
Suppose $\P$ is a poset and $A\subseteq\P$. The {\bf upper shadow} of $A$ in $\P$ is
\[
    \uSdw_\P A = \{p\in\P \mid p\text{ covers }a\text{ for some }a\in A\}
\] and the {\bf lower shadow} of $A$ in $\P$ is
\[
    \dSdw_\P A = \{p\in\P \mid a\text{ covers }p\text{ for some }a\in A\}\,.
\]
\end{defn}
\begin{ex}\label{ex: upper shadow} Consider the poset of monomials in the polynomial ring $K[x,y]$ over a field $K$ with partial order given by divisibility. We identify a monomial $x^ay^b$ with the unit square with lower left corner at $(a,b)$ in the cartesian plane. In the figure below, under this identification, the set $A=\{x^3,y^3\}$ is shown in pink and its upper shadow $\uSdw A=\{x^4,x^3y,y^3x, y^4\}$ is shown in green.
\begin{figure}[h!]
    \centering
    \includegraphics[height=4cm]{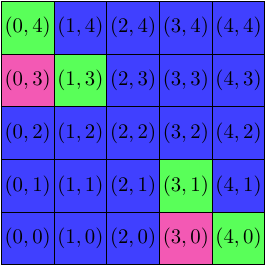}
    \caption{The monomial poset of $K[x,y]$ with $\{x^3,y^3\}$ and its upper shadow highlighted.}
\end{figure}
\end{ex}

 A {\bf rank function} on a poset $\P$ with a minimum element $m$ is a function $r: \P\rightarrow\mathbb{N}$ such that $r(p)+1=r(q)$ whenever $q$ covers $p$, and $r(m)=0$. A poset $\P$ is {\bf ranked} if and only if there exists a rank function on $\P$. If a poset is ranked, then its rank function is unique. If $d\in\mathbb{N}$ and $r$ is a rank function on $\P$, then let $\P_d = r^{-1}(d)$.

\begin{defn}
Suppose $\P$ is a ranked poset. Suppose $\prec$ is a total order on $\P$. If $d,s\in\mathbb{N}$ with $s\leq\lvert\P_d\rvert$, then the {\bf initial segment} of size $s$ in $\P_d$ is the set consisting of the largest $s$ elements of $\P_d$ with respect to $\prec$, and is denoted $\Seg_d s$.
\end{defn}

\begin{defn}
A ranked poset $\P$ is {\bf Macaulay} if and only if there exists a total order  on $\P$ such that for all subsets $A\subseteq\P_d$, the following hold
\begin{enumerate}
\item Initial segments have the smallest upper shadows: 
\[
\left | \uSdw_\P\Seg_d |A|  \right | \leq  |\uSdw_\P(A)|;
\]
\item The upper shadow of an initial segment is an initial segment:
\[
 \uSdw_\P  \Seg_d |A| =\Seg_{d+1} |\uSdw_\P\Seg_d|A||.
\]
\end{enumerate}
\end{defn}

\begin{ex}\label{ex: upper shadows} We consider the poset \(\cP\) of monomials in \(K[x,y]\), ordered by divisibility, with the lexicographic order \(\O\) as a total order. On the left, the set \(A = \{x^3, xy^2\}\) (pink) and its upper shadow (green) are shown, yielding \(|A| = 2\) and \(|\uSdw(A)| = 4\). On the right, the set \(B =\Seg_2\lvert A\rvert = \{x^3, x^2y\}\) consists of the two lexicographically largest monomials of degree two, with its upper shadow (green) satisfying \(|\uSdw(B)| = 3 \leq |\uSdw(A)| = 4\). Notably, \(\uSdw(B)\) itself forms a segment under \(\O\), containing the three largest lexicographic monomials of degree three in \(\cP\).

\begin{figure}[h!]
\centering
\includegraphics{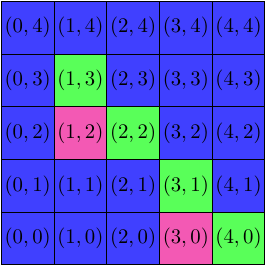}
\qquad
\includegraphics{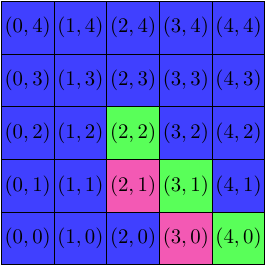}
\caption{Two sets (pink) and their upper shadows (green) in the monomial poset of $K[x,y]$.}
\label{fig: upper shadows}
\end{figure}
\end{ex}
\subsection{Macaulay rings}\label{s: MacaulayRings}

Let $n\in\mathbb{N}$, and let $R = K[x_1,\dots,x_n]$ be a polynomial ring in $n$ variables over a field $K$. Suppose $I\subseteq R$ is a homogeneous ideal.  A {\bf monomial} in $R/I$ is a nonzero  coset $m+I$ for some monomial $m$ in $R$.

\begin{defn}\label{def: monomial poset}
Let $S$ be a graded quotient ring of a polynomial ring and let $\M_S$ be the set of monomials in $S$. Define a partial order denoted $\vert$ on $\M_S$ by $f\vert g$ if and only if $fh=g$ for some $h\in\M_S$. The pair $(\M_S, \vert)$ is termed the {\bf monomial poset} of $S$. 
\end{defn}

An example of a monomial poset is presented in \Cref{fig: poset}.

%The ring $S$ is {\bf level linearly independent} if and only if for each $d\in\mathbb{N}$, the set of degree $d$ monomials in $S$ is linearly independent over $K$.

%Suppose $S$ is level linearly independent and $J\subseteq S$ is an ideal. 
The {\bf Hilbert function} of a homogeneous ideal $I$ of $S$ is the function $H_I: \mathbb{N}\rightarrow\mathbb{N}$ defined by  $H_I(d)=\dim_K I_d$, i.e. the dimension  of the vector space of degree $d$ elements of $I$.

\begin{defn}
A graded quotient ring $S$ of a polynomial ring over a field is {\bf Macaulay} if and only if there exists a total order $\O$ on $\M_S$ such that for each homogeneous ideal $I\subset S$ the vector space $\O^*(I)$ described below is an ideal of $S$
\[
   \O^*(I)= \bigoplus_{d\geq 0}\Span_K\Seg_dH_I(d)
\]
and additionally $H_I(d)=H_{\O^*(I)}(d)$ for all $d\in \N$.
\end{defn}

It is shown in \cite[Theorem 2.6.3]{Kuz} that a ring is Macaulay if and only if its monomial poset is Macaulay. Macaulay's theorem \cite{Mac} says that polynomial rings $K[x_1, \dots, x_n]$ are Macaulay for all $n\in\N$. The Clements-Lindstr\"om Theorem \cite{CL} says that monomial complete intersection rings, i.e., rings of the form $\frac{K[x_1, \dots, x_n]}{(x_1^{d_1}, \dots, x_n^{d_n})}$ for $d_i\in\N$, are Macaulay. Few other encompassing classes of Macaulay rings are known. For some further examples the reader is invited to peruse \cite{BL, Kuz, Polymath paper}.

%\subsection{Software for studying posets}

%\red{We should say that there is already a Posets package, cite the paper about it \cite{PosetsJSAG}, and mention it does not have functions that check if a poset is Macaulay or deal with monomial posets from non-monomial quotient rings. We should also mention Bezrukovs Java tool.}

%The package \texttt{Posets} \cite{PosetsJSAG} provides a data type and many functions for working with posets. It has a function \texttt{standardMonomialPoset} which can compute the monomial poset of a monomial ideal, but it has is no function for computing posets for non-monomial ideals. It also has no function for checking whether a poset is Macaulay.

%There is a software tool \cite{SoftwareTool} which can check whether a poset with at most $30$ elements is Macaulay. 

\section{Using the MacaulayPosets package}

First, we load the package. We can also set \texttt{printWidth} to a sufficiently small number if we want to ensure long outputs do not spill into the margins of this paper.

\begin{verbatim}
    i1 : loadPackage "MacaulayPosets"
    i2 : printWidth = 70
\end{verbatim}

\subsection{Monomial posets}

One feature of interest in the package is the ability to construct monomial posets from commutative rings using the command \texttt{getPoset}. As discussed in \cref{s: MacaulayRings}, our notion of divisibility in graded rings is quite general. In particular it applies to quotients of a polynomial ring by ideals which are not necessarily monomial. For monomial ideals, the command \texttt{standardMonomialPoset} offered by the \texttt{Posets} package returns the monomial poset of the corresponding quotient ring, but cannot handle non-monomial ideals. Our main contribution presented in this section is to develop code which is applicable to any graded ring that computes its monomial poset in the sense of \Cref{def: monomial poset}.

\begin{ex}
Let $R = \Q[x, y]$ be a polynomial ring, $I = ( x^3, y^5)$, a monomial ideal in $R$. One can use the function \verb|getMons| to obtain the monomials of $R/I$:
\begin{verbatim}
    i3 : R = QQ[x, y]
    i4 : I = ideal(x^3, y^5)
    i5 : getMons(R, I)
                 2   2    2 2   2 3   2 4          2     3     4      2
    o5 : {1, x, x , x y, x y , x y , x y , x*y, x*y , x*y , x*y , y, y ,
         ----------------------------------------------------------------
          3   4
         y , y }
\end{verbatim}

Note that the last input could be substituted by any of the following: 
\begin{compactitem}
    \item \verb|getMons(R/I)|
    \item \verb|getMons(I)| (given that Macaulay2 is using $R$),
    \item \verb|getMons(QQ[x, y]/ideal(x^3, y^5))|. 
\end{compactitem}
\end{ex}

\begin{ex}
Let $R = \Q[x, y]$ be a polynomial ring, $I = (xy - y^2, x^4, x^3y, x^2y^2)$, a homogeneous ideal that is not a monomial ideal, and $S = R/I$. 
To find the set $\M_S$, we use the code: 
\begin{verbatim}
    i3 : getMons( QQ[x, y]/ideal(x*y - y^2, x^4, x^3 * y, x^2 * y^2  ) )
                 2   3      2   3
    o4 : {1, x, x , x , y, y , y }
\end{verbatim}
\end{ex}

There are instances where we have quotient rings that have infinitely many monomials. We can use the optional input \texttt{MaxDegree} to limit to what degree what monomials we want. The default value is 10. The optional input \texttt{MaxDegree} will be ignored when the quotient ring entered in has a finite number of monomials. 

\begin{ex}
Let us consider the quotient ring $S = \mathbb{Q}[x, y]/(x^2 - y^2)$ which has infinitely many monomials.
\begin{verbatim}
    i3 : getMons( QQ[x, y], ideal(x^2 - y^2), MaxDegree => 8 )
                    2          2   3   4     3     4   5   6     5     6
    o4 : {1, x, y, y , x*y, x*y , y , y , x*y , x*y , y , y , x*y , x*y ,
         ----------------------------------------------------------------
          7   8     7
         y , y , x*y }
\end{verbatim}
\end{ex}

Next we illustrate the usage of the \verb|getPoset| function to find the monomial poset of various graded rings under the divisibility partial order. This function also has the optional input \texttt{MaxDegree} with analogous use as in \texttt{getMons}.

\begin{ex}
Let $R = \Q[x, y], I = ( x^2, y^2 ), S = R/I$. 

The poset of monomials can be gotten from the following code:
\begin{verbatim}
    i3 : QQ[x,y]
    i4 : I = monomialIdeal(x^2, y^2)
    i5 : P = getPoset I
\end{verbatim}
This poset can also be obtained by using \texttt{standardMonomialPoset} from the \texttt{Posets} package.
\begin{verbatim}
    i6 : Q = standardMonomialPoset I
    i7 : areIsomorphic( P, Q )
    o7 : true
\end{verbatim} 
\end{ex}

The \texttt{Posets} package has no equivalent method function for ideals which are not monomial. However the function \texttt{getPoset} in the package \texttt{MacaulayPosets} does not exhibit such restrictions.

\begin{ex}
Let $R = \Q[x, y], I = (x^6, x^3y, y^4, x^2y^3 - x^5), S = R/I$ \\
We can get the poset $\M_S$ by the following code: 
\begin{verbatim}
    i3 : getPoset( QQ[x, y], ideal(x^6, x^3 * y, y^4, x^2*y^3 - x^5) )
\end{verbatim}

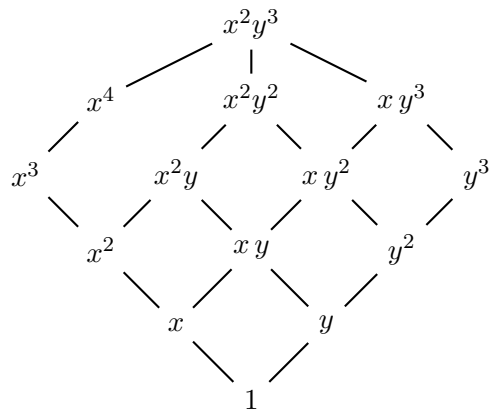
\begin{figure}[h!]
\centering
\begin{tikzpicture}[xscale=2, every path/.style={thick}]
        draw \node at (0, 0) (0) {$1$};
        draw \node at (-1/2, 1) (1) {$x$};
        draw \node at (1/2, 1) (2) {$y$};
        draw \node at (-1, 2) (3) {$x^{2}$};
        draw \node at (0, 2) (4) {$x\,y$};
        draw \node at (1, 2) (5) {$y^{2}$};
        draw \node at (-3/2, 3) (6) {$x^{3}$};
        draw \node at (-1/2, 3) (7) {$x^{2}y$};
        draw \node at (1/2, 3) (8) {$x\,y^{2}$};
        draw \node at (3/2, 3) (9) {$y^{3}$};
        draw \node at (-1, 4) (10) {$x^{4}$};
        draw \node at (0, 4) (11) {$x^{2}y^{2}$};
        draw \node at (1, 4) (12) {$x\,y^{3}$};
        draw \node at (0, 5) (13) {$x^{2}y^{3}$};

        \draw (0) -- (1);
        \draw (0) -- (2);
        \draw (1) -- (4);
        \draw (1) -- (3);
        \draw (3) -- (7);
        \draw (3) -- (6);
        \draw (6) -- (10);
        \draw (10) -- (13);
        \draw (7) -- (11);
        \draw (11) -- (13);
        \draw (4) -- (7);
        \draw (4) -- (8);
        \draw (8) -- (11);
        \draw (8) -- (12);
        \draw (12) -- (13);
        \draw (2) -- (5);
        \draw (2) -- (4);
        \draw (5) -- (9);
        \draw (5) -- (8);
        \draw (9) -- (12);
\end{tikzpicture}

\caption{The monomial poset of $\mathbb{Q}[x,y]/(x^6,x^3y,y^4,x^2y^3-x^5)$.}
\label{fig: poset}
\end{figure}
\end{ex}

\subsection{Lower and Upper Shadows}
The functions \texttt{lowerShadow} and \texttt{upperShadow} will compute the respective sets as described in \Cref{def: shadows}. 

\begin{ex} One can compute \Cref{ex: upper shadow} as follows.
\begin{verbatim}
    i3 : R = QQ[x,y]
    i4 : x = R_0
    i5 : y = R_1
    i6 : upperShadow(getPoset(R, MaxDegree=>4), {x^3, y^3})
           4   3      3   4
    o6 : {x , x y, x*y , y }
\end{verbatim}

\end{ex}

\subsection{Additivity}
For this section, there is only one function to mention: \texttt{isAdditive}. 
Additivity was introduced in \cite{ClementsAdditive} and was leveraged in \cite{Clements} to study disjoint unions of isomorphic posets. Important classes of posets, such as the monomial posets of polynomial rings and monomial complete intersection rings, are additive. 
Before we introduce this property, we need the following definitions.

%\begin{defn}
		Suppose that $\mathcal{P}$ is a ranked poset with an additional total order $\prec$ on it.
		A set $A$ consisting of elements of the same rank is called a {\bf segment}, 
		if for any $a,b\in A$ and any $c\in \mathcal{P}$ of the same rank such that $a \prec c \prec b$,
		we must have $c\in A$. The {\bf initial segment} of size $q$ in the $d$-th level $\P_d$ of $\P$ consists of  the largest $q$ elements with respect to $\prec$ in $P_d$, while the {\bf final segment} of size $q$ in $\P_d$ is  the smallest $q$ elements with respect to $\prec$ in $P_d$.  
%\end{defn}

\begin{ex}
	Suppose that we consider the poset of monomials of $K[x,y]$ with the lex order.
	The set $ \{x^4y^4,x^5y^3,x^6y^2\}$ is a segment, 
	but $\{x^4y^4,x^6y^2\}$ is not a segment.
\end{ex}

A finite sequence of integers is called additive if the sum of any consecutive number of its terms is less than or equal to the sum of the same number of initial terms in the sequence and greater than or equal to the sum of the same number of final terms in the sequence. Additivity for posets, as defined below, is equivalent to the sequence  $\{|\{\uSdw_{\text{new}} (p_i)\}|\}_i$ being additive where $p_1\succ \cdots \succ p_s$ is the ordered sequence of elements of $\P_d$ for a fixed $d$. It is shown in \cite[p.40]{ClementsAdditive} that additive posets satisfy subadditivity of upper shadows of initial segments. Specifically for all positive integers $a,b,d$ the following inequality holds
\[
|\uSdw(\Seg_d(a+b))|\leq |\uSdw(\Seg_d a)| + |\uSdw(\Seg_d b)|.
\]

\begin{defn}
		Suppose that $\mathcal{P}$ is a ranked poset with an additional total order $\prec$ on it.
		Let $A$ be a segment of some rank and $B$ be all the elements of the same rank that are larger than all the elements of $A$.
		The {\bf new shadow} of $A$ is defined to be
		\begin{align*}
			\uSdw_{\text{new}} (A) = \uSdw(A)\setminus \uSdw(B).
		\end{align*}
	\end{defn}
\begin{defn}\label{def: additive}
		Suppose that $\mathcal{P}$ is a Macaulay poset.
		We say that $\mathcal{P}$ is {\bf additive} if the following hold:
		\begin{enumerate}
			\item If $A$ is an initial segment and $B$ is a segment such that $|A|=|B|$ then 
			$$|\uSdw_{\text{new}}(A)| \geq |\uSdw_{\text{new}}(B)|.$$
			\item If $B$ is a segment and $C$ is a final segment such that $|B|=|C|$ then 
			$$|\uSdw_{\text{new}}(B)| \geq |\uSdw_{\text{new}}(C)|.$$
		\end{enumerate}
	\end{defn}
\begin{ex}

The new shadow of the initial segment $A = \{x^3, x^2y\}$ from \Cref{ex: upper shadows} (where it is called $B$) is shown in \Cref{fig: new shadows} on the left. Its upper shadow is shown in \Cref{fig: upper shadows} on the right. The rank $3$ monomials larger than all elements of $A$ are $\{xy^2, y^3\}$. So, the new shadow of $A$ is
\begin{align*}
    \uSdw_{\text{new}} (A) &= \uSdw A \setminus \uSdw\{xy^2, y^3\} \\
    &= \{x^4, x^3y, x^2y^2\}\setminus \{x^2y^2, xy^3, y^4\}\\
    &= \{x^4, x^3y\}\,.
\end{align*}
The new shadows of the segment $B = \{x^2y, xy^2\}$ and the final segment $C = \{xy^2, y^3\}$ are shown in the middle and right of \Cref{fig: new shadows}, respectively. For the three sets $A,B,C$, the inequalities \[\lvert \uSdw_{\text{new}}(A) \rvert \leq \lvert \uSdw_{\text{new}}(B)\rvert \leq \lvert \uSdw_{\text{new}}(C) \rvert\]
from \Cref{def: additive} hold, i.e. $2\leq 2\leq 3$ holds. The monomial poset of any polynomial ring, in particular $K[x,y]$, is known to be additive.

\begin{figure}[h!]
\centering
\includegraphics{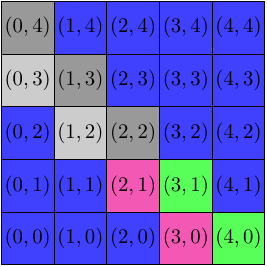}
\quad
\includegraphics{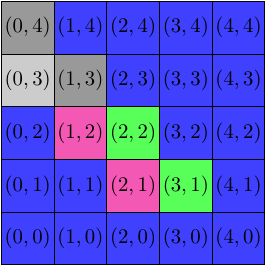}
\quad
\includegraphics{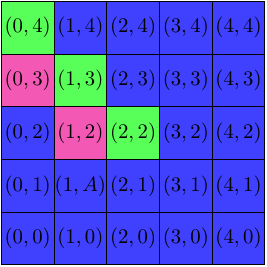}
\caption{Three lexicographic segments (pink) and their new shadows (green) in the monomial poset of $K[x,y]$. The elements in light gray are larger than all those in pink and their upper shadow is marked in dark gray. }
\label{fig: new shadows}
\end{figure}
\end{ex}
\begin{ex} The function \texttt{isAdditive} detects this property.
\begin{verbatim}
    i3 : isAdditive product(chain 4, chain 4)
    o3 : true
    i4 : isAdditive getPoset( QQ[x,y,z]/(x^4,y^2-z^2,z^2-x*y) )
    o4 : false
\end{verbatim}

\end{ex}

\subsection{Macaulay posets and rings}

The command \texttt{isMacaulay} can be used to determine whether a poset or ring is Macaulay.

\begin{ex}\label{ex: 39 vertices}
Below we verify that the $2^5$-element Boolean lattice is Macaulay but that its disjoint union with the $7$-element chain is not. By \cite[Theorem 2.6.3]{Kuz}, another way to verify that this Boolean lattice is Macaulay is by applying \texttt{isMacaulay} to the ring $\frac{\mathbb{Q}[v,w,x,y,z]}{(v^2,w^2,x^2,y^2,z^2)}$. The option \texttt{Visual=>true} will open an interactive browser window, where TikZ code compiling to \Cref{fig: booleanLattice union chain} was obtained.

\begin{verbatim}
    i3 : isMacaulay booleanLattice 5
    o3 : true
    i4 : isMacaulay(union(booleanLattice 5, chain 7), Visual=>true)
    o4 : false
    i5 : isMacaulay( QQ[v,w,x,y,z]/(v^2,w^2,x^2,y^2,z^2) )
    o5 : true
\end{verbatim}

\begin{figure}[h!]
    \scalebox{0.75}{\begin{tikzpicture}
    \newcommand*\pointsZYdgy{285.3333333333333/564/0/00000,122.28571428571429/475/1/00001,244.57142857142858/475/2/00010,71.33333333333333/386/3/00011,366.8571428571429/475/4/00100,142.66666666666666/386/5/00101,214/386/6/00110,71.33333333333333/297/7/00111,489.14285714285717/475/8/01000,285.3333333333333/386/9/01001,356.66666666666663/386/10/01010,142.66666666666666/297/11/01011,428/386/12/01100,214/297/13/01101,285.3333333333333/297/14/01110,122.28571428571429/208/15/01111,611.4285714285714/475/16/10000,499.3333333333333/386/17/10001,570.6666666666666/386/18/10010,356.66666666666663/297/19/10011,642/386/20/10100,428/297/21/10101,499.3333333333333/297/22/10110,244.57142857142858/208/23/10111,713.3333333333333/386/24/11000,570.6666666666666/297/25/11001,642/297/26/11010,366.8571428571429/208/27/11011,713.3333333333333/297/28/11100,489.14285714285717/208/29/11101,611.4285714285714/208/30/11110,285.3333333333333/119/31/11111,570.6666666666666/564/32/1,733.7142857142858/475/33/2,784.6666666666666/386/34/3,784.6666666666666/297/35/4,733.7142857142858/208/36/5,570.6666666666666/119/37/6,428/30/38/7}
    \newcommand*\edgesZYdgy{0/1,0/2,0/4,0/8,0/16,1/3,1/5,1/9,1/17,2/3,2/6,2/10,2/18,3/7,3/11,3/19,4/5,4/6,4/12,4/20,5/7,5/13,5/21,6/7,6/14,6/22,7/15,7/23,8/9,8/10,8/12,8/24,9/11,9/13,9/25,10/11,10/14,10/26,11/15,11/27,12/13,12/14,12/28,13/15,13/29,14/15,14/30,15/31,16/17,16/18,16/20,16/24,17/19,17/21,17/25,18/19,18/22,18/26,19/23,19/27,20/21,20/22,20/28,21/23,21/29,22/23,22/30,23/31,24/25,24/26,24/28,25/27,25/29,26/27,26/30,27/31,28/29,28/30,29/31,30/31,32/33,33/34,34/35,35/36,36/37,37/38}
    \newcommand*\scaleZYdgy{0.02}
    \foreach \x/\y/\z/\w in \pointsZYdgy
      \node (\z) at (\scaleZYdgy*\x,-\scaleZYdgy*\y) [circle,draw,inner sep=0pt] {$\w$};
    \foreach \x/\y in \edgesZYdgy
      \draw (\x) -- (\y);
    \end{tikzpicture}}

\caption{The TikZ picture generated by the \texttt{Visualize} package for the non-Macaulay disjoint union of the $2^5$-element Boolean lattice with the $7$-element chain.}
\label{fig: booleanLattice union chain}
\end{figure}
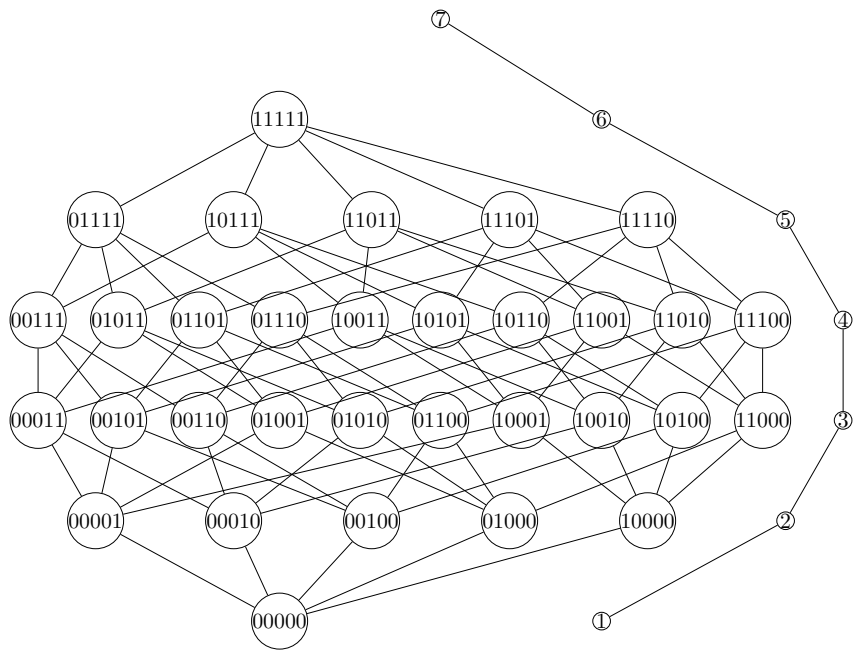

\end{ex}
  The function \texttt{isMacaulay} includes two optional arguments: \texttt{TikZ} and \texttt{Visual}.
\begin{compactitem}
\item If the \texttt{TikZ} option is set to \texttt{true}, the function generates TikZ code, allowing the user to visualize the poset of monomials in \LaTeX. The vertices are horizontally arranged according to a Macaulay order.
\item
Setting the \texttt{Visual} option to \texttt{true} enables direct visualization of the poset and its Macaulay Orders within Macaulay2. This relies on the \texttt{Visualize} package \cite{Visualize}.
\end{compactitem}

The function \texttt{macaulayOrders} returns a list of total orders with respect to which the given poset  is Macaulay. This function has options \texttt{TikZ}, \texttt{Visual} and \texttt{AllOrders}. The first two work the same way as with \texttt{isMacaulay}. With \texttt{true} as the default value, \texttt{AllOrders} will give a list of all possible total orders of the poset that are Macaulay.

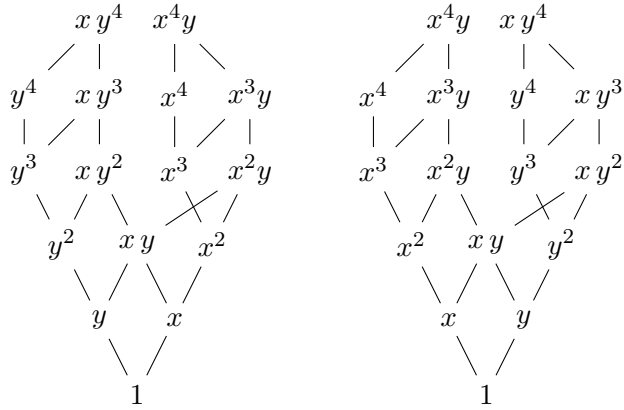
\begin{figure}[ht!]
\centering 
\begin{tikzpicture}
        draw \node at (0, 0) (0) {$1$};
        draw \node at (-1/2, 1) (1) {$y$};
        draw \node at (1/2, 1) (2) {$x$};
        draw \node at (-1, 2) (3) {$y^{2}$};
        draw \node at (0, 2) (4) {$x\,y$};
        draw \node at (1, 2) (5) {$x^{2}$};
        draw \node at (-3/2, 3) (6) {$y^{3}$};
        draw \node at (-1/2, 3) (7) {$x\,y^{2}$};
        draw \node at (1/2, 3) (8) {$x^{3}$};
        draw \node at (3/2, 3) (9) {$x^{2}y$};
        draw \node at (-3/2, 4) (10) {$y^{4}$};
        draw \node at (-1/2, 4) (11) {$x\,y^{3}$};
        draw \node at (1/2, 4) (12) {$x^{4}$};
        draw \node at (3/2, 4) (13) {$x^{3}y$};
        draw \node at (-1/2, 5) (14) {$x\,y^{4}$};
        draw \node at (1/2, 5) (15) {$x^{4}y$};

        \draw (0) -- (1);
        \draw (0) -- (2);
        \draw (2) -- (4);
        \draw (2) -- (5);
        \draw (5) -- (9);
        \draw (5) -- (8);
        \draw (8) -- (12);
        \draw (8) -- (13);
        \draw (12) -- (15);
        \draw (13) -- (15);
        \draw (9) -- (13);
        \draw (4) -- (7);
        \draw (4) -- (9);
        \draw (7) -- (11);
        \draw (11) -- (14);
        \draw (1) -- (4);
        \draw (1) -- (3);
        \draw (3) -- (7);
        \draw (3) -- (6);
        \draw (6) -- (11);
        \draw (6) -- (10);
        \draw (10) -- (14);
\end{tikzpicture}
\qquad
\begin{tikzpicture}
        draw \node at (0, 0) (0) {$1$};
        draw \node at (-1/2, 1) (1) {$x$};
        draw \node at (1/2, 1) (2) {$y$};
        draw \node at (-1, 2) (3) {$x^{2}$};
        draw \node at (0, 2) (4) {$x\,y$};
        draw \node at (1, 2) (5) {$y^{2}$};
        draw \node at (-3/2, 3) (6) {$x^{3}$};
        draw \node at (-1/2, 3) (7) {$x^{2}y$};
        draw \node at (1/2, 3) (8) {$y^{3}$};
        draw \node at (3/2, 3) (9) {$x\,y^{2}$};
        draw \node at (-3/2, 4) (10) {$x^{4}$};
        draw \node at (-1/2, 4) (11) {$x^{3}y$};
        draw \node at (1/2, 4) (12) {$y^{4}$};
        draw \node at (3/2, 4) (13) {$x\,y^{3}$};
        draw \node at (-1/2, 5) (14) {$x^{4}y$};
        draw \node at (1/2, 5) (15) {$x\,y^{4}$};

        \draw (0) -- (2);
        \draw (0) -- (1);
        \draw (1) -- (4);
        \draw (1) -- (3);
        \draw (3) -- (7);
        \draw (3) -- (6);
        \draw (6) -- (10);
        \draw (6) -- (11);
        \draw (10) -- (14);
        \draw (11) -- (14);
        \draw (7) -- (11);
        \draw (4) -- (9);
        \draw (4) -- (7);
        \draw (9) -- (13);
        \draw (13) -- (15);
        \draw (2) -- (4);
        \draw (2) -- (5);
        \draw (5) -- (9);
        \draw (5) -- (8);
        \draw (8) -- (13);
        \draw (8) -- (12);
        \draw (12) -- (15);
\end{tikzpicture}

\caption{The two orders with respect to which the monomial poset of $\Q[x, y]/(x^5, x^2y^2, y^5)$ from \Cref{ex: macaulayOrders} is Macaulay, with smaller elements with respect to the total order appearing to the left of larger elements.}
\label{fig: macaulayOrders}
\end{figure}

\begin{ex}\label{ex: macaulayOrders}
Consider the monomial poset of $\Q[x, y]/(x^5, x^2y^2, y^5)$.
One find its Macaulay total orders in the following way: 
\begin{verbatim}
    i3 : macaulayOrders(QQ[x,y]/(x^5,x^2*y^2,y^5), TikZ=>true)
    o4 : {{1, y, x, y , x*y, x , y , x*y , x , x y, y , x*y , x , x y,
         ----------------------------------------------------------------
            4   4               2        2   3   2    3     2   4   3
         x*y , x y}, {1, x, y, x , x*y, y , x , x y, y , x*y , x , x y,
         ----------------------------------------------------------------
          4     3   4      4
         y , x*y , x y, x*y }}
\end{verbatim}
The printed TikZ code compiles to \Cref{fig: macaulayOrders}. This shows an example of a heart-shaped poset equipped with the twist order, as discussed in section 5 of \cite{Polymath paper}.

Using\texttt{macaulayOrders(P, Visual=>true)} one can get an interactive Hasse diagram of the poset $P$ and all Macaulay orders on $P$. The function \texttt{macaulayOrders(P, Visual=>true)} will open a new tab in the default browser with visualization options for the poset $P$. If the end session button is clicked, then another tab with similar visualization options will open with the first Macaulay order on $P$. If the end session button is clicked on the second tab that was opened, then a third tab will open with the next Macaulay order on $P$. We can continue like this until we have cycled through all the Macaulay orders on $P$. If no new tab is opened after we click end session for the first time then it means that there are no Macaulay orders on $P$.
\end{ex}

\subsection{Operations on posets}\label{s: poset ops} Four poset operations are included in \texttt{MacaulayPosets}. The first three are discussed in detail in \cite{Polymath paper}. 

\begin{defn}
Suppose that for $1 \leq i \leq t$ we have posets $P_i$ each with unique least element $\ell_i$. Their \textbf{wedge product} is the set: 
\[ P_1 \vee P_2 \vee \cdots \vee P_t  = \left(\bigsqcup_{i=1}^t P_i \right)/ (\ell_1=\ell_2=\cdots  =\ell_t),  \] meaning that we take the disjoint union of the sets $P_i$ in which we identify all the $\ell_i$ into one element, with the partial order $a\leq b$ if and only if $a\leq b$ in $P_i$ for some $i$.
\end{defn}
\begin{figure}[h!]
\centering 
\begin{tikzpicture}[scale=1, every path/.style={thick}, every node/.style={circle, draw, inner sep=2.5pt}]
	draw \node at (0, 0) (0) {};
	draw \node at (-1, 1) (1) {};
	draw \node at (0, 1) (2) {};
	draw \node at (1, 1) (3) {};
	draw \node at (-1, 2) (4) {};
	draw \node at (0, 2) (5) {};
    draw \node at (1, 2) (6) {};
	draw \node at (-1, 3) (7) {};
	draw \node at (0, 3) (8) {};
    draw \node at (1, 3) (9) {};

	\draw (0) -- (1) -- (4) -- (7);
	\draw (0) -- (2) -- (5) -- (8);
	\draw (0) -- (3) -- (6) -- (9);
\end{tikzpicture}
\caption{The wedge product of three $4$-element chains.}
\label{fig: wedge}
\end{figure}
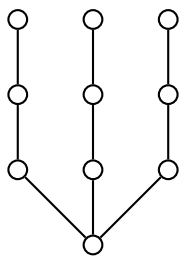

To calculate the wedge product of several posets, all of which have a unique least element, we use the function: \verb|posetWedgeProduct|. 

\begin{ex}
Let $P$ be the $4$-element chain. To obtain the poset $P\lor P\lor P$, shown in \Cref{fig: wedge} one would use the \verb|chain| function in the \verb|Posets| package to get path posets of varying lengths and then enter the following code: 
\begin{verbatim}
    i3 : posetWedgeProduct({chain 4, chain 4, chain 4})
\end{verbatim}
\end{ex}

\begin{ex}
We can also place posets obtained from wedge products into a wedge product: 
\begin{verbatim}
    i3 : P1 = posetWedgeProduct( chain(4), booleanLattice(3) )
    i4 : P2 = posetWedgeProduct( getPoset( QQ[x, y]/ideal(x^3, y^4) ),
        getPoset( QQ[a, b]/ideal( a^3, a*b, b^2 ) ) )
    i5 : posetWedgeProduct({P1, P2, chain(7), booleanLattice(2)})
\end{verbatim}
\end{ex}

Wedge products are the special case of the following operation with $P$ a singleton poset.

\begin{defn}
Suppose $P$ is a poset and for $1 \leq i \leq t$ we have posets $Q_i$ with a rank-preseving injective monotone map $f_i: P\rightarrow Q_i$. The \textbf{fiber product} of the $Q_i$ over $P$ is the set: 
\[ Q_1 \times_P Q_2 \times_P \cdots \times_P Q_t  = \left(\bigsqcup_{i=1}^t Q_i \right)/ (f_i(p)=f_j(p)),  \] meaning that we take the disjoint union of the sets $Q_i$ in which, for each $p\in P$, we identify all the $f_i(p)$ into one element, with the partial order $a\leq b$ if and only if $a\leq b$ in $Q_i$ for some $i$.
\end{defn}

Fiber products of posets can be constructed using \texttt{posetFiberProduct}. Its parameters are \texttt{PosetMap}s, which can be obtained from \texttt{map} by supplying a codomain, a domain, and images of the vertices of the domain.

\begin{ex}
The monomial poset shown in \Cref{fig: macaulayOrders} is isomorphic to a fiber product.

\begin{verbatim}
    i3 : P = product(chain 2, chain 2)
    i4 : Q = product(chain 2, chain 5)
    i5 : f = map(Q, P, {{1,1},{1,2},{2,1},{2,2}})
    i6 : g = map(Q, P, {{1,1},{2,1},{1,2},{2,2}})
    i7 : areIsomorphic( posetFiberProduct(f,g), getPoset(QQ[x,y]/(x^5,x^
    2*y^2,y^5)) )
    o7 : true
\end{verbatim}
\end{ex}

Another operation implemented in this package is the closed product of a list of posets. In the paper \cite{Polymath paper} we term this operation diamond product to match the symbol used to denote it. However, the \texttt{Posets} package implements a different operation under the name diamond, so we opted for a terminology that avoids potential confusion.

\begin{defn}
Suppose that for $1 \leq i \leq t$ we have posets $P_i$ each with unique least element $\ell_i$ and unique maximum element $L_i$. Their \textbf{closed product} is the set: 
\[ P_1 \diamond P_2 \diamond \cdots \diamond  P_t = \left(\bigsqcup_{i=1}^t P_i \right)/ (\ell_1=\ell_2=\cdots \ell_t, L_1=L_2=\cdots=L_t ), \] meaning that we take the disjoint union of the sets $P_i$ in which we identify all the $\ell_i$ into one element and all the $L_i$ elements into one element, with the partial order $a\leq b$ if and only if $a\leq b$ in $P_i$ for some $i$.
\end{defn}
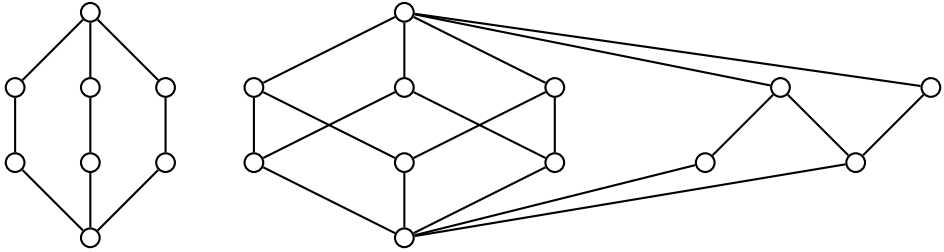
\begin{figure}[h!]
\centering 
\begin{tikzpicture}[scale=1, every path/.style={thick}, every node/.style={circle, draw, inner sep=2.5pt}]
	draw \node at (0, 0) (0) {};
	draw \node at (-1, 1) (1) {};
	draw \node at (0, 1) (2) {};
	draw \node at (1, 1) (3) {};
	draw \node at (-1, 2) (4) {};
	draw \node at (0, 2) (5) {};
    draw \node at (1, 2) (6) {};
	draw \node at (0, 3) (7) {};

	\draw (0) -- (1) -- (4) -- (7);
	\draw (0) -- (2) -- (5) -- (7);
	\draw (0) -- (3) -- (6) -- (7);
\end{tikzpicture}
\qquad
\begin{tikzpicture}[scale=1, every path/.style={thick}, every node/.style={circle, draw, inner sep=2.5pt}]
    draw \node at (2,0) (000) {};
       
    draw \node at (0,1) (001) {};
    draw \node at (2,1) (010) {};
    draw \node at (4,1) (100) {};
        
    draw \node at (0,2) (011) {};
    draw \node at (2,2) (101) {};
    draw \node at (4,2) (110) {};

    draw \node at (2,3) (111) {};

    draw \node at (6, 1) (w) {};
    draw \node at (8, 1) (x) {};

    draw \node at (7, 2) (y) {};
    draw \node at (9, 2) (z) {};

    \draw (000)--(001)--(011)--(111)--(101);
    \draw (010)--(000)--(100)--(110)--(111);
    \draw (001)--(101)--(100);
    \draw (011)--(010)--(110);

    \draw (000)--(w)--(y)--(111);
    \draw (000)--(x)--(z)--(111);
    \draw (x)--(y);
\end{tikzpicture}
\caption{Examples of closed products.}
%\caption{The closed product of three $4$-element chains.}
\label{fig: closed}
\end{figure}

The closed products  shown in \Cref{fig: closed} are obtained in \Cref{ex: closed prod} below.

\begin{ex}\label{ex: closed prod}
\hspace{0em}
\begin{verbatim}
    i3 : posetClosedProduct({chain 4, chain 4, chain 4})
    i4 : posetClosedProduct( booleanLattice(3), getPoset( QQ[x, y]/ideal(x^2, y^3)) )
\end{verbatim}

In case the posets used are not of the same rank, the resulting poset will not be ranked. An accompanying error message will be printed.

\begin{verbatim}
    i5 : posetClosedProduct({booleanLattice(2), getPoset( QQ[x, y]/ideal(x^6,
    y^3) ), chain(5)})
    The given posets do not have the same rank, therefore the resulting poset
    will not be ranked.
\end{verbatim}
\end{ex}

Closed products are the special case of the following operation with $P$ a singleton poset.

\begin{defn}
Suppose $P$ is a poset and for $1 \leq i \leq t$ we have self-dual posets $Q_i$ with a rank-preseving injective monotone map $f_i: P\rightarrow Q_i$. Let $q_i: Q_i\rightarrow Q_i^{\operatorname{op}}$ be the isomorphisms of the $Q_i$ with their duals. The \textbf{connected sum} of the $Q_i$ over $P$ is the set: 
\[ Q_1 \#_P Q_2 \#_P \cdots \#_P Q_t  = \left(\bigsqcup_{i=1}^t Q_i \right)/ (f_i(p)=f_j(p), (q_i\circ f_i)(p) = (q_j\circ f_j)(p)),  \] meaning that we take the disjoint union of the sets $Q_i$ in which, for each $p\in P$, we identify all the $f_i(p)$ into one element and all the $(q_i\circ f_i)(p)$ into one element, with the partial order $a\leq b$ if and only if $a\leq b$ in $Q_i$ for some $i$.
\end{defn}

\begin{ex}
Here is the connected sum of a $2\times 5$ box with itself over a $2$-element chain.
\begin{center}
\begin{tikzpicture}[scale=0.9]

\node at (-8,0) {\begin{tikzpicture}[xscale=-2, scale=0.7, every path/.style={thick}, every node/.style={circle, draw, inner sep=2.5pt}]
        draw \node[blue] at (1/2, 0) (0) {};
        draw \node[blue] at (0, 1) (2) {};
        draw \node at (1, 1) (3) {};
        draw \node at (1/2, 2) (6) {};
        draw \node at (3/2, 2) (7) {};
        draw \node[red] at (2, 3) (9) {};
        draw \node at (1, 3) (10) {};
        draw \node[red] at (3/2, 4) (11) {};

        \draw[blue] (0) -- (2);
        \draw (0) -- (3);
        \draw (3) -- (7);
        \draw (3) -- (6);
        \draw (6) -- (10);
        \draw (2) -- (6);
        \draw (7) -- (9);
        \draw[red] (9) -- (11);
        \draw (7) -- (10);
        \draw (10) -- (11);
\end{tikzpicture}
};

\node at (-5.5,0) {\scalebox{1.2}{$\#$}};

\node at (-3,0) {\begin{tikzpicture}[xscale=2, scale=0.7, every path/.style={thick}, every node/.style={circle, draw, inner sep=2.5pt}]
        draw \node[blue] at (1/2, 0) (0) {};
        draw \node[blue] at (0, 1) (2) {};
        draw \node at (1, 1) (3) {};
        draw \node at (1/2, 2) (6) {};
        draw \node at (3/2, 2) (7) {};
        draw \node[red] at (2, 3) (9) {};
        draw \node at (1, 3) (10) {};
        draw \node[red] at (3/2, 4) (11) {};

        \draw[blue] (0) -- (2);
        \draw (0) -- (3);
        \draw (3) -- (7);
        \draw (3) -- (6);
        \draw (6) -- (10);
        \draw (2) -- (6);
        \draw (7) -- (9);
        \draw[red] (9) -- (11);
        \draw (7) -- (10);
        \draw (10) -- (11);
\end{tikzpicture}
};

\node at (-0.5,0) {\scalebox{1.2}{$=$}};

\node at (3,0) {
\begin{tikzpicture}[xscale=2, scale=0.7, every path/.style={thick}, every node/.style={circle, draw, ,inner sep=2.5pt}]
        draw \node[blue] at (0, 0) (0) {};
        draw \node at (-1, 1) (1) {};
        draw \node[blue] at (0, 1) (2) {};
        draw \node at (1, 1) (3) {};
        draw \node at (-3/2, 2) (4) {};
        draw \node at (-1/2, 2) (5) {};
        draw \node at (1/2, 2) (6) {};
        draw \node at (3/2, 2) (7) {};
        draw \node at (-1, 3) (8) {};
        draw \node[red] at (0, 3) (9) {};
        draw \node at (1, 3) (10) {};
        draw \node[red] at (0, 4) (11) {};

        \draw (0) -- (1);
        \draw[blue] (0) -- (2);
        \draw (0) -- (3);
        \draw (3) -- (7);
        \draw (3) -- (6);
        \draw (6) -- (10);
        \draw[line width=2pt, white] (7) -- (9);
        \draw (7) -- (9);
        \draw[red] (9) -- (11);
        \draw (2) -- (6);
        \draw (2) -- (5);
        \draw (7) -- (10);
        \draw (10) -- (11);
        \draw (1) -- (4);
        \draw (1) -- (5);
        \draw (5) -- (8);
        \draw[line width=2pt, white] (4) -- (9);
        \draw (4) -- (9);
        \draw (4) -- (8);
        \draw (8) -- (11);
\end{tikzpicture}
};

\end{tikzpicture}
\end{center}
This can be obtained using \texttt{posetConnectedSum}.

\begin{verbatim}
    i3 : P = chain 2
    i4 : Q = product(chain 2, chain 4)
    i5 : f = map(Q, P, {{1,1}, {2,1}})
    i6 : posetConnectedSum(f,f)
    \end{verbatim}

\end{ex}

\subsection{Operations on rings}

Two ring operations are included in \texttt{MacaulayPosets}, namely the fiber product and connected sum. When the monomial posets of the resulting rings are considered, these operations are the algebraic counterpart of the poset operations termed wedge product, poset fiber product,  closed product, and poset connected sum in the previous section.

\begin{defn}\label{def: ring FP}
Given rings $A, B$, and $C$ equipped with ring homomorphisms $\pi_A:A\to C$ and $\pi_B:B\to C$, the {\bf fiber product} of $A$ and $B$ over $C$ is the following subring of $A\times B$ 
\[
A\times_C B=\{(a,b): a\in A, b\in B, \pi_A(a)=\pi_B(b)\}.
\]

In the particular case when $C=K$, $A = R/I$  and $B=S/J$ for some homogeneous ideals $I$ of $R = K[x_1,\ldots, x_n]$ and $J$ of $S=K[y_1, \ldots, y_m]$, the {\bf fiber product over $K$} admits an explicit presentation
\begin{equation}\label{eq: presentation FP}
A\times_K B= \frac{K[x_1,\ldots, x_n, y_1, \ldots, y_m]}{ I +J + (x_iy_j: 1\leq i\leq n, 1\leq j\leq m)}.
\end{equation}
\end{defn}

Since a presentations for the general fiber product are difficult to write down explicitly, in our package we offer the function \texttt{ringFiberProduct} which computes the fiber product of two  $K$-algebras over $K$ using the presentation \eqref{eq: presentation FP}.

\begin{ex}
Consider the quotient rings: 
\[
A = \Q[x, y, z]/( x^4, x^2 y^2, z^3 ) \text{ and } B = \Q[a, b]/(a^6, a^2 b^2, b^3).
\]
We can compute the ring fiber product over $\Q$ using the following code: 
\begin{verbatim}
    i3 : ringFiberProduct( QQ[x, y, z]/ideal(x^4, x^2 * y^2, z^3),
    QQ[a, b]/(a^6, a^2 * b^2, b^3) )
     
                                QQ[x..z, a..b]
    o3 = ----------------------------------------------------------
            4   2 2   3   6   2 2   3
          (x , x y , z , a , a b , b , x*a, x*b, y*a, y*b, z*a, z*b)
\end{verbatim}
\end{ex}

\begin{defn}\label{def: ring CS}
Given Gorenstein rings $A, B$  and $C$ equipped with ring homomorphisms $\pi_A:A\to C$ and $\pi_B:B\to C$ and dual maps $\iota_A:C\to A$, $\iota_B:C\to B$ such that $(\iota_A(c),\iota_B(c))\in A\times_C B$ for all $c\in C$, the {\bf connected sum} of $A$ and $B$ over $C$ is the ring $A\#_C B$ defined by means of the following exact sequence \cite[(2.3)]{IMS}
\[
0\longrightarrow C \xrightarrow{(\iota_A, \iota_B)} A\times_C B\longrightarrow A\#_C B\longrightarrow 0.
\]

Suppose now that $A = R/I$  and $B=S/J$ for some homogeneous ideals $I$ of $R = K[x_1,\ldots, x_n]$ and $J$ of $S=K[y_1, \ldots, y_m]$ and $C=K$. Then $A$ and $B$ have unique (up to scalar) maximal elements $m_A$ and $m_B$ respectively, called {\em socle elements}. The  {\bf connected sum} of $A$ and $B$ over $K$ admits a presentation
\begin{equation}\label{eq: presentation CS}
A\#_K B= \frac{K[x_1,\ldots, x_n, y_1, \ldots, y_m]}{I +J + (x_iy_j: 1\leq i\leq n, 1\leq j\leq m)+(m_A-m_B)}.
\end{equation}
\end{defn}

The \texttt{MacaulayPosets} package offer the function \texttt{ringConnectedSum} which computes the connected sum of two  $K$-algebras over $K$ using the presentation \eqref{eq: presentation CS}.
Our code does not check whether the rings input are Gorenstein or not. We instead check only whether the monomial posets of each quotient ring have unique maximal elements or not and pick the elements $m_A$ and $m_B$ in \eqref{eq: presentation CS} to be the maximal elements of the posets $\M_A$ and $\M_B$ with respect to the divisibility partial order, respectively. 

\begin{ex} 
Let $A = \Q[x, y]/(x^3, y^4)$ and let $B = \Q[a, b]/(a^4, b^2)$. The connected sum $A \#_{\Q} B$ is obtained with the code: 
\begin{verbatim}
    i3 : ringConnectedSum(QQ[x, y]/ideal(x^3, y^4 ), QQ[a, b]/ideal(a^4, b^2))

                          QQ[x..y, a..b]
    o3 = ------------------------------------------------
           3   4   4   2                       2 3    3
         (x , y , a , b , x*a, x*b, y*a, y*b, x y  - a b)
\end{verbatim}
\end{ex}

The ring and poset operations discussed in this paper are inspired by  topological  constructions. In topology, the fiber product  of topological spaces combines two spaces over a common base space based on given maps between them and the connected sum of two spaces is a new space formed by removing an open neighborhood  from each space and gluing the resulting boundaries together. The ring-theoretic constructions considered in this section are abstractions of these ideas. Indeed, in the category of rings  the fiber product in \Cref{def: ring FP} and in the category of Gorenstein rings the connected sum in \Cref{def: ring CS} applied to  cohomology rings of two spaces being combined by one of the previously mentioned topological operations yield the cohomology ring of the fiber product and respectively of the connected sum spaces. 

Passing from graded rings to their monomial posets  yields the poset operations described in the previous sections of the paper. The precise correspondence between the ring operations described above and the poset operations in \cref{s: poset ops} 
can be summarized as follows: the monomial poset of a fiber product of graded rings is the fiber product of the monomial posets of the summands and the monomial poset of a connected sum of graded Gorenstein rings  is the connected sum of the monomial posets of the summands. Specializing to fiber products and connected sums over the field $K$ gives a correspondence between these operations and the poset operations of wedge product and closed product as follows:
\begin{eqnarray*}
\M_{A\times_K B}=\M_A \vee\M_B  \\
\M_{A\#_K B}=\M_A \diamond\M_B.
\end{eqnarray*}
The first displayed equation 
is proved in \cite[Proposition 2.22]{Polymath paper} and the second follows by similar considerations.

\end{document}